\DeclareFontFamily{OT1}{rsfs}{}
\DeclareFontShape{OT1}{rsfs}{n}{it}{<-> rsfs10}{}
\DeclareMathAlphabet{\mathscr}{OT1}{rsfs}{n}{it}
\documentclass[11pt,amssymb]{amsart}
\usepackage{amsmath}
\usepackage{amssymb}
\usepackage{mathabx}
\usepackage[colorlinks]{hyperref}
\def\cD{{\mathscr D}}
\def\cE{{\mathscr E}}

\def\ocM{\overline{\mathscr M}}
\def\oM{\overline M}

\def\cG{\mathscr G}

\def\cR{{\mathscr R}}

\def\cT{{\mathscr T}}

\def\bQ{{\mathbb Q}}

\def\fd{{\mathfrak d}}
\def\fD{{\mathfrak D}}
\def\ff{{\mathfrak f}}
\def\fF{{\mathfrak F}}

\def\cT{{\mathscr T}}

\def\cC{{\mathscr C}}
\begin{document}
\newtheorem {theo}{Theorem}
\newtheorem {coro}{Corollary}
\newtheorem {lemm}{Lemma}
\newtheorem {rem}{Remark}
\newtheorem {defi}{Definition}
\newtheorem {ques}{Question}
\newtheorem {prop}{Proposition}
\def\spb{\smallpagebreak}
\def\mpb{\vskip 0.5truecm}
\def\bpb{\vskip 1truecm}
\def\wtM{\widetilde M}
\def\tM{\widetilde M}
\def\wtN{\widetilde N}
\def\tN{\widetilde N}
\def\tC{\widetilde C}
\def\tX{\widetilde X}
\def\tY{\widetilde Y}
\def\ti{\widetilde \iota}
\def\bs{\bigskip}
\def\ms{\medskip}
\def\ni{\noindent}
\def\td{\nabla}
\def\pd{\partial}
\def\hol{$\text{hol}\,$}
\def\Log{\mbox{Log}}
\def\bfQ{{\bf Q}}
\def\Todd{\mbox{Todd}}
\def\bP{{\bf P}}
\def\dxi{d x^i}
\def\dxj{d x^j}
\def\dyi{d y^i}
\def\dyj{d y^j}
\def\dzi{d z^I}
\def\dzj{d z^J}
\def\ozi{d{\overline z}^I}
\def\ozj{d{\overline z}^J}
\def\oz1{d{\overline z}^1}
\def\oz2{d{\overline z}^2}
\def\oz3{d{\overline z}^3}
\def\sI{\sqrt{-1}}
\def\hol{$\text{hol}\,$}
\def\ok{\overline k}
\def\ol{\overline l}
\def\oJ{\overline J}
\def\oT{\overline T}
\def\oS{\overline S}
\def\oV{\overline V}
\def\oW{\overline W}
\def\oI{\overline I}
\def\oK{\overline K}
\def\oL{\overline L}
\def\oj{\overline j}
\def\oi{\overline i}
\def\ok{\overline k}
\def\oz{\overline z}
\def\om{\overline mu}
\def\on{\overline nu}
\def\oa{\overline \alpha}
\def\ob{\overline \beta}
\def\of{\overline f}
\def\og{\overline \gamma}
\def\ogamma{\overline \gamma}
\def\oGamma{\overline\Gamma}
\def\odelta{\overline \delta}
\def\otheta{\overline \theta}
\def\ophi{\overline \phi}
\def\opd{\overline \partial}
\def\oA{\overline A} 
\def\oB{\overline B}
\def\oC{\overline C}
\def\op{\overline D}
\def\oIq1{\oI_1\cdots\oI_{q-1}}
\def\oIq2{\oI_1\cdots\oI_{q-2}}
\def\op{\overline \partial}
\def\ua{{\underline {a}}}
\def\us{{\underline {\sigma}}}
\def\tor{{\mbox{tor}}}
\def\vol{{\mbox{vol}}}
\def\rank{{\mbox{rank}}}
\def\bp{{\bf p}}
\def\bk{{\bf k}}
\def\a{{\alpha}}
\def\tchi{\widetilde{\chi}}


\title[Arithmetic fake compact Hermitian symmetric spaces of Type $A_3$] 
{Arithmetic fake compact Hermitian symmetric spaces of Type $A_3$}
 \maketitle
{\centerline{\sc Gopal Prasad and Sai-Kee Yeung}}

\bs

\begin{center}
{\bf 1. Introduction}
\end{center}
\vskip3mm

This paper is a sequel to our paper [PY3] and we expect that the reader is familiar with the basic definitions and techniques of the earlier paper. The notion of {\it arithmetic fake compact Hermitian symmetric space} has been defined in 1.2 of [PY3]. 
It has been shown in [PY2] and [PY3] that arithmetic fake compact Hermitian symmetric spaces can only be of type $A_1$, $A_2$, $A_3$ and $A_4$. The ones of type $A_2$ are the {\it fake projective planes} which are smooth projective complex 
algebraic surfaces of considerable interest. These have been classified in [PY1] (please see the corrected version of this paper posted on arXiv) in $28$ nonempty classes, and we know now that up to biholomorphism there are $100$ fake projective planes. It has been shown in [PY2] that there are some arithmetic fake compact Hermitian symmetric spaces of type $A_4$, for example a few arithmetic fake ${\bf{P}}^4_{\bf C}$ and also a few arithmetic fake ${\bf{Gr}}_{2,5}$.  Using the techniques of [PY1] and the results of [PY2] it should not be hard to determine all arithmetic fake compact Hermitian symmetric spaces of type $A_4$. 
On the other hand, there are too many examples of such spaces of type $A_1$ for a convenient classification. Thus the only arithmetic fake compact Hermitian symmetric spaces which possibly can be classified, but have remained to be classified,  
are the ones of type $A_3$. The computations for $A_3$ is considerably more difficult than the ones handled in [PY3].  The goal of this article is to investigate spaces of type $A_3$.  To this end, we show that there are none  except possibly those arising from the pairs of number fields  listed in Theorem 1 at the end of Section 6.  In these possible cases, one may hope to complete the classification in the same way as the classification of fake projective planes achieved in [PY1], [CS].  However, at this moment it appears that the computational power of available softwares is not sufficient to carry out the required computation.   In sections 7 and 8, we list some questions that may help to decipher the possibilities in the future.   
\vskip2mm

\ni{\bf 1.1}  We begin by giving a description of groups of type $A_3$ of interest to us in this paper. Let $k$ be a number field of degree $d$ over $\bQ$ and $V_f$ be the set of nonarchimedean places of $k$. Let $G$ be an absolutely almost simple  simply connected $k$-group of type $A_3$. Let $v_1,\,\ldots,\,v_r$ be the archimedean places of $k$ where $G$ is isotropic (equivalently, $G(k_{v_i})$ is noncompact) and let  $\mathscr{G} = \prod_{i=1}^{r}G(k_{v_i})$ considered as a real Lie group with with the product topology. Let $X$ denote the symmetric space of $\cG$ and $X_u$ be the compact dual of $X$. We assume that $r\geqslant 1$ and the symmetric space $X$ is Hermitian. Then for $i\leqslant r$, $G(k_{v_i})$ is isomorphic to ${\rm{SU}}(m_i, 4-m_i)$, with $m_i\ne 0$, hence each $v_i$ is real. Moreover, for any archimedean place $v$ of $k$ where $G$ is anisotropic, $G(k_v)$ is isomorphic to the compact group ${\rm{SU}}(4)$, so  any such $v$ is also real. Thus we see that $k$ is totally real and $G$ is of type $^{2}A_3$, i.e., it is an outer form of a split group. Let $\ell$ be the quadratic extension of $k$ over which $G$ is inner.  Then $\ell$ is totally complex. 
\vskip1mm

From the classification of classical groups we know that there is a division algebra $\cD$ with center $\ell$ and of degree $\mathfrak{d}$ dividing 4, $\cD$ given with an involution $\sigma$ of second kind such that $k$ is fixed pointwise under $\sigma$, and a hermitian form $h$ on $\cD^{4/{\mathfrak{d}}}$ defined in terms of the involution $\sigma$  so that $G$ is the special unitary group ${\rm{SU}}(h)$.   Since the division algebra $\cD$  is equipped with an involution $\sigma$ of second kind, 
if $v$ is a place of $k$ which does not split in $\ell$, then 
$k_v\otimes_k \cD$ is the matrix algebra $M_{\mathfrak d}(\ell_v)$, where $\ell_v$ denotes the field $\ell\otimes_k k_v$. 
\vskip1mm

If $\fd\ne 1$, i.e., $\cD \ne \ell$,  then $\fd = 4$ or $2$. Let $v$ be a place of $k$ which splits in $\ell$, then there is a division algebra $\fD_v$ with center $k_v$ and of degree $\fd_v$, $\fd_v\,\vert\, \fd$, such that $k_v\otimes_k \cD = M_{\fd/\fd_v}(\fD_v)\times M_{\fd/\fd_v}(\fD_v^{o})$, where $\fD_v^{o}$ is the opposite of $\fD_v$, and the involution $\sigma$ interchanges the two factors of $k_v\otimes_k \cD$.  For such a $v$, $G(k_v)$ is isomorphic to ${\rm{SL}}_{\fd/\fd_v}(\fD_v)$. 
\vskip1mm

Let $\cT_0$ be the set of places $v$ of $k$ which split in $\ell$ and $\fd_v>1$. From Class Field Theory we know that only for finitely many $v$, $\fd_v>1$, and moreover, there exists $v$ such that $\fd_v = \fd$, so  $\cT_0$ is finite and nonempty (if $\fd \ne 1$).  As $k$ is totally real and $\ell$ is totally complex, none of the archimedean places of $k$ split  in $\ell$. So every $v\in \cT_0$ is nonarchimedean. 
\vskip1mm

\ni {\bf 1.2} Let $\overline{G}$ be the adjoint group of $G$ and let $\pi: G\rightarrow \overline{G}$ be the natural isogeny. The kernel of this isogeny is described in 1.5 of [PY3] for groups of type $^2A_n$; for $n=3$, $C({\bf R})$ is a cyclic group of order $4$.  We will denote the image $\pi(\cG)\subseteq \prod_{i=1}^r{\overline G}(k_{v_i})$ by $\overline{\cG}$. The subgroup $\overline\cG$ is the identity component of $ \prod_{i=1}^r{\overline G}(k_{v_i})$. Let $\Pi\subset {\overline{\cG}}$ be the fundamental group of an arithmetic fake  compact Hermitian symmetric space (arithmetic with respect to the $k$-group structure on $G$) which is a compact quotient of $X$. Then $\Pi$ is a torsion-free co-compact discrete subgroup of $\overline{\cG}$, and $\Pi$ is arithmetic with respect to the $k$-group structure on $G$. Let $\widetilde{\Pi}$ be the inverse image of $\Pi$ in  $\cG$. Then  as the kernel of $\widetilde{\Pi}\rightarrow \Pi$ is a subgroup of order $4^r$, and the Euler-Poincar\'e characteristic $\chi(\Pi)$ of $\Pi$ equals $\chi(X_u)$, the Euler-Poincar\'e characteristic $\chi({\widetilde{\Pi}})$ of $\widetilde{\Pi}$\,(in the sense of C.\,T.\,C.\,Wall) equals $\chi(X_u)/4^r$. Now let $\Gamma$ be a maximal discrete subgroup of $\cG$ containing $\widetilde{\Pi}$ and $\Lambda = \Gamma\cap G(k)$, $G(k)$ embedded in $\cG=\prod_{i=1}^r G(k_{v_i})$ diagonally. The subgroup $\Lambda$ is a ``principal'' arithmetic subgroup, i.e., for every nonarchimedean place $v$ of $k$, the closure $P_v$ of $\Lambda$ in $G(k_v)$ is a parahoric subgroup. By the strong approximation property of the simply connected semi-simple group $G$ we see that $\Lambda = G(k)\cap \prod P_v$, moreover $\Gamma$ is the normalizer of $\Lambda$ in $\cG$; see Proposition 1.4 of [BP] for these results.   In terms of the Haar measure $\mu$ used in [P]  and [BP], and to be used here, $\chi({\widetilde{\Pi}})=\chi(X_u) \mu(\cG/{\widetilde{\Pi}})$. So we conclude that $\mu(\cG/{\widetilde{\Pi}})= 1/4^r$ and therefore $\mu(\cG/\Gamma)$ is a submultiple of $1/4^r$. Hence, in particular,  $4^r\mu(\cG/\Gamma) \leqslant 1$.  
\vskip1mm

\ni {\bf 1.3} Let $\cT$ be the set of nonarchimedean places $v$ of $k$ such that {\it either} (i) $v$ does not ramify in $\ell$ (equivalently, $G$ splits over an unramified extension of $k_v$) and $P_v$ is not a hyperspecial parahoric subgroup of $G(k_v)$, {\it or} (ii) $v$ ramifies in $\ell$ and either $G$ is not quasi-split over $k_v$ (i.e., its $k_v$-rank is $1$) or $P_v$ is not special.  (Note that $\cT$ is precisely the set of nonarchimedean places $v$ of $k$ such that $e'(P_v)\ne 1 $. It contains the $\cT$  introduced in \S2.10 of [PY3].) As every place  $v\in \cT_0$ splits in $\ell$, but $G$ does not split over $k_v$ (and hence $G(k_v)$ cannot contain a hyperspecial parahoric subgroup) we see that $\cT_0\subseteq \cT$.  

\begin{center}
{\bf 2. Preliminaries}
\end{center}
\vskip3mm

\ni{\bf 2.1}  All unexplained notation are from [P], [BP] and [PY3]. The value of 
$\mathfrak{s}(\mathscr{G})$ given in [P] for $^2A_{3}$ is 
$5.$ 
Hence the formula for the covolume of the principal arithmetic group $\Lambda= G(k)\cap \prod P_v$ given in 
[P] is: 

\begin{equation}
\mu(\cG/\Lambda )= (D_kD_\ell)^{5/2}\Big(\frac3{2^7\pi^9}\Big)^{d}\cE, 
\end{equation}
where
\begin{eqnarray}
\cE&=&\prod_{v\in V_f}e(P_v)\\
&=&\zeta_k(2)L_{\ell|k}(3)\zeta_k(4)\prod_{v\in V_f} e^{\prime}(P_v). 
\end{eqnarray}

The index of $\Lambda$ in the maximal arithmetic $\Gamma$ is bounded 
by (see (54) in \S8 of [PY3])
\begin{equation}
[\Gamma:\Lambda]\leqslant 2^{d+2r+2\#\cT}h_{\ell, 4},
\end{equation}
where $h_{\ell,4}$ is the order of the subgroup 
of the class group of $\ell$ consisting of elements of order dividing $4$.

\bs
\ni{\bf 2.2} 
The values  of $e(P_v)$ and $e'(P_v)$ appearing in (2) and (3) are
$$e(P_v) = \frac{q_v^{(\dim \oM_v +\dim \ocM_v)/2}}{\#\oM_v(\ff_v)},$$
$$e'(P_v)=e(P_v)\cdot\frac{\#\ocM_v(\ff_v)}{q_v^{\dim\ocM_v}},$$
in the notation of  [P].  
In the following we first list the possible values of $\frac{\#\ocM_v(\ff_v)}{q_v^{\dim\ocM_v}}$.  

\ms
\ni{\it Case 1:} $v$ splits in $\ell$.  Then $\ocM_v=SL_4$ and 
$$\frac{\#\ocM_v(\ff_v)}{q_v^{\dim\ocM_v}}=(1-\frac1{q_v^2})(1-\frac1{q_v^3})(1-\frac1{q_v^4}).$$

\ms
\ni{\it Case 2:} $v$ is inert in $\ell$, i.e., it does not split in $\ell$, but is unramified.  Then $\ocM_v=SU_4$ and 
$$\frac{\#\ocM_v(\ff_v)}{q_v^{\dim\ocM_v}}=(1-\frac1{q_v^2})(1+\frac1{q_v^3})(1-\frac1{q_v^4}).$$

\ms
\ni{\it Case 3:} $v$ is ramified in $\ell$.  Then  $\ocM_v$ is of type $C_2$ and 
$$\frac{\#\ocM_v(\ff_v)}{q_v^{\dim\ocM_v}}=(1-\frac1{q_v^2})(1-\frac1{q_v^4}).$$

\bs
\ni{\bf 2.3} Now we list the values of $e(P_v)$ and $e'(P_v)$. 

\vskip2mm
\ni{\it Case 1:} {\it $v$ splits in $\ell$ and $G$ splits over $k_v$.  There are five possibilities:}
\vskip1mm

{\it Case 1a:} $P_v$ is hyperspecial, then $\oM_v= {\rm{SL}}_4$, and in this case 
$$e(P_v)=(1-\frac1{q_v^2})^{-1}(1-\frac1{q_v^3})^{-1}(1-\frac1{q_v^4})^{-1}, \ e'(P_v)=1.$$

{\it Case 1b:} $\oM_v$ is isogenous to ${\rm{GL}}_3$, in which case
$$e(P_v)=q_v^3(1-\frac1{q_v})^{-1}(1-\frac1{q_v^2})^{-1}(1-\frac1{q_v^3})^{-1}, \  e'(P_v)=(q_v+1)(q_v^2+1).$$

{\it Case 1c:} $\oM_v$ is isogenous to ${\rm{GL}}_1\times ({\rm{SL}}_2)^2$, in which case 
$$e(P_v)=q_v^4(1-\frac1{q_v})^{-1}(1-\frac1{q_v^2})^{-2}, \  e'(P_v)=(q_v^2+1)(q_v^2+q_v+1).$$

{\it Case 1d:} $\oM_v$ is isogenous to $({\rm{GL}}_1)^2\times {\rm{SL}}_2$, in which case 
$$e(P_v) = q_v^5(1-\frac 1{q_v^2})^{-1}(1-\frac 1{q_v})^{-2}, \ e'(P_v) = (q_v+1)(q_v^2+1)(q_v^2+q_v+1).$$

{\it Case 1e:} $P_v$ is a Iwahori subgroup, then  $\oM_v=({\rm{GL}}_1)^3$ and 
$$e(P_v)=q_v^6(1-\frac1{q_v})^{-3}, \ e'(P_v)=(q_v+1)^2(q_v^2+1)(q_v^2+q_v+1).$$

\vskip2mm
\ni{\it Case 2:} {\it $v$ splits in $\ell$ and $G$ does not split over $k_v$.  There are three possibilities:}
\vskip1mm

{\it Case 2a:} $\rank_{k_v}G=0$.  In this case, $\oM_v=R_{\fF_v/\ff_v}({\rm{GL}}_1)/{\rm{GL}}_1$, where $\fF_v$ is the extension of degree $4$ of $\ff_v$.
In this case,
$$e(P_v)=q_v^6(1+\frac1{q_v})^{-1}(1+\frac1{q_v^2})^{-1}, \  e'(P_v)=(q_v-1)(q^2_v-1)(q^3_v-1).$$

{\it Case 2b(i):} $\rank_{k_v}G=1$ and $P_v$ is a maximal parahoric subgroup. In this case $\oM_v$ is isogenous to the product of $R_{\fF_v/\ff_v}({\rm{SL}}_2)$ with the unique 1-dimensional $\ff_v$-anisotropic torus $R^{(1)}_{\fF_v/\ff_v}({\rm{GL}}_1)$, where $\fF_v$ is the  extension of degree $2$ of $\ff_v$, and 
$$e(P_v)=q_v^4(1+\frac1{q_v})^{-1}(1-\frac1{q_v^4})^{-1}, \  e'(P_v)=(q_v-1)(q_v^3-1).$$

{\it Case 2b(ii):} $\rank_{k_v}G =1$ and $P_v$ is an Iwahori subgroup. In this case $\oM_v$ is isomorphic to $R_{\fF_v/\ff_v}(({\rm{GL}}_1)^2)/{\rm{GL}}_1$, where $\fF_v$ is the extension of degree $2$ of $\ff_v$ and 
$$e(P_v)=q_v^6(1+\frac1{q_v})^{-1}(1-\frac1{q_v^2})^{-1}, \  e'(P_v)=(q_v-1)(q^2_v+1)(q^3_v-1).$$

\vskip2mm
\ni{\it Case 3:} {\it $v$ is inert in $\ell$.  There are two possibilities:}
\vskip1mm

{\it Case 3a:} $\rank_{k_v}G=2$ (i.e.\,$G$ is quasi-split over $k_v$).  There are five possible subcases.
\vskip1mm

\hspace{0.4cm} {\it Case 3a(i):} $P_v$ is hyperspecial, in which case $\oM_v=\ocM_v$ and these groups are isomorphic to ${\rm{SU}}_4$ over $\ff_v$.  In this case,
$$e(P_v)=(1-\frac1{q_v^2})^{-1}(1+\frac1{q_v^3})^{-1}(1-\frac1{q_v^4})^{-1}, \ e'(P_v)=1.$$

\hspace{0.4cm}{\it Case 3a(ii):} $\oM_v$ is isgenous  to the product of  ${\rm{SL}}_2\times {\rm{SL}}_2$ with the 1-dimensional $\ff_v$-anisotropic torus.  In this case 
$$e(P_v)=q_v^4(1+\frac1{q_v})^{-1}(1-\frac1{q_v^2})^{-2}, \ e'(P_v)=(q_v^2+1)(q_v^2-q_v+1).$$

\hspace{0.4cm}{\it Case 3a(iii):} 
$\oM_v$ is isogenous to the product $R_{\fF_v/\ff_v}({\rm{SL}}_2)\times {\rm{GL}}_1$ in which case 
$$e(P_v)=q_v^4(1-\frac1{q_v})^{-1}(1-\frac1{q_v^4})^{-1}, \  e'(P_v)=(q_v+1)(q^3_v+1).$$

\hspace{0.4cm}{\it Case 3a(iv):} $\oM_v$ is isogenous to the product of ${\rm{SL}}_2$ and the torus $R_{\fF_v/\ff_v}({\rm{GL}}_1)$, where $\fF_v$ is the extension of degree 2 of $\ff_v$.  In this case,
$$e(P_v)=q_v^5(1-\frac1{q_v^2})^{-2}, \ e'(P_v)=(q_v^2+1)(q_v^3+1).$$

\hspace{0.4cm}{\it Case 3a(v):} $P_v$ is an Iwahori subgroup. Then $\oM_v$ is isogenous to the product $R_{\fF_v/\ff_v}({\rm{GL}}_1)\times {\rm{GL}}_1$, where $\fF_v$ is the degree 2 extension of $\ff_v$.  In this case,
$$e(P_v)=q_v^6(1-\frac1{q_v})^{-1}(1-\frac1{q_v^2})^{-1}, \  e'(P_v)=(q_v+1)(q^2_v+1)(q^3_v+1).$$

{\it Case 3b:} $\rank_{k_v}G=1$.  There are two subcases.
\vskip1mm

\hspace{0.4cm}{\it Case 3b(i):} $P_v$ is a maximal parahoric subgroup. In this case $\oM_v$ is isogenous to the product of ${\rm{SU}}_3$ with the 1-dimensional $\ff_v$-anisotropic torus, and 
$$e(P_v)=q_v^3(1+\frac1{q_v})^{-1}(1-\frac1{q_v^2})^{-1}(1+\frac1{q_v^3})^{-1}, \  e'(P_v)=(q_v-1)(q^2_v+1).$$

\hspace{0.4cm}{\it Case 3b(ii):} $P_v$ is an Iwahori subgroup. Then $\oM_v=R_{\fF_v/\ff_v}({\rm{GL}}_1)^2/{\rm{GL}}_1$, where $\fF_v$ is the degree 2 extension of $\ff_v$.  In this case,
$$e(P_v)=q_v^6(1+\frac1{q_v})^{-1}(1-\frac1{q_v^2})^{-1}, \  e'(P_v)=(q_v-1)(q^2_v+1)(q^3_v+1).$$

\vskip2mm
\ni{\it Case 4:} {\it $v$ is ramified in $\ell$.  There are two possibilities:}
\vskip1mm

{\it Case 4a:} $\rank_{k_v}G=2$ (i.e. $G$ is quasi-split over $k_v$).  There are four subcases.

\hspace{0.4cm}{\it Case 4a(i):} $P_v$ is a special maximal parahoric subgroup. Then $\oM_v=\ocM_v$ and these groups are isomorphic to ${\rm{Sp}}_4$ (i.e., are of type $C_2$).  In this case,
$$e(P_v)=(1-\frac1{q_v^2})^{-1}(1-\frac1{q_v^4})^{-1}, \ e'(P_v)=1.$$

\hspace{0.4cm}{\it Case 4a(ii):} $\oM_v$ is isogenous to ${\rm{SL}}_2\times {\rm{SL}}_2$.  In this case,
$$e(P_v)=q_v^2(1-\frac1{q_v^2})^{-2}, \  e'(P_v)=(q^2_v+1).$$

\hspace{0.4cm}{\it Case 4a(iii):} $\oM_v$ is isogenous to ${\rm{GL}}_2$.  In this case,
$$e(P_v)=q_v^3(1-\frac1{q_v})^{-1}(1-\frac1{q_v^2})^{-1}, \  e'(P_v)=(q_v+1)(q^2_v+1).$$

\hspace{0.4cm}{\it Case 4a(iv):} $P_v$ is an Iwahori subgroup. Then $\oM_v={\rm{GL}}_1\times {\rm{GL}}_1$.  In this case,
$$e(P_v)=q_v^4(1-\frac1{q_v})^{-2}, \  e'(P_v)=(q_v+1)^2(q_v^2+1).$$

{\it Case 4b:} $\rank_{k_v}G=1$.  There are three subcases.
\vskip1mm

\hspace{0.4cm}{\it Case 4b(i):} $\oM_v=R_{\fF_v/\ff_v}({\rm{SL}}_2)$, where $\fF_v$ is the degree 2 extension of $\ff_v$.  In this case,
$$e(P_v)=q_v^2(1-\frac1{q_v^4})^{-1}, \  e'(P_v)=(q^2_v-1).$$

\hspace{0.4cm}{\it Case 4b(ii):} $\oM_v$ is isogenous to the product of ${\rm{SL}}_2$ and the 1-dimensional $\ff_v$-anisotropic torus.  In this case,
$$e(P_v)=q_v^3(1+\frac1{q_v})^{-1}(1-\frac1{q_v^2})^{-1}, \  e'(P_v)=(q_v-1)(q^2_v+1).$$

\hspace{0.4cm}{\it Case 4b(iii):} $P_v$ is an Iwahori subgroup. Then $\oM_v=R_{\fF_v/\ff_v}({\rm{GL}}_1)$, where $\fF_v$ is the degree 2 extension of $\ff_v$.  In this case,
$$e(P_v)=q_v^4(1-\frac1{q_v^2})^{-1}, \  e'(P_v)=(q^4_v-1).$$

\bs
\ni{\bf 2.4}
From the functional equations $$\zeta_k(2j) = D_k^{\frac{1}{2}-
2j}\big(\frac{(-1)^{j}2^{2j-1}\pi^{2j}}{(2j-1)!}\big)^d \zeta_k(1-2j),$$
and $$L_{\ell|k}(2j+1) = \big(\frac{D_k}{D_{\ell}}\big)^{2j+\frac{1}{2}}\big(\frac{(-1)^j 2^{2j}\pi^{2j+1}}{(2j)!}\big)^d L_{\ell|k}(-2j),$$
we find that 
\begin{equation}
\cR :=(D_k D_\ell)^{5/2}\Big(\frac3{2^7\pi^9}\Big)^{d}\big(
\zeta_k(2)L_{\ell|k}(3)\zeta_k(4)\big)
\end{equation}
$$\ \ \ \ =2^{-3d}\zeta_k(-1)L_{\ell|k}(-2)\zeta_k(-3).$$
For all nonarchimedean $v\notin \cT$, as $e'(P_v) =1$, quations (1), (2) and (3) give that  
\begin{equation*}
\mu(\mathcal{G}/\Lambda)= \cR\prod_{v\in \cT}e'(P_v).
\end{equation*}
As the values of $e'(P_v)$ given in {1.3} are integral, we conclude that $\mu(\cG/\Lambda)$ 
is an integral multiple of $\cR$.  Moreover,

\begin{equation}
\mu(\cG/\Gamma)= \frac{\mu(\cG/\Lambda)}{[\Gamma :\Lambda ]}= \frac{\cR \prod_{v\in \cT}e'(P_v)}{[\Gamma :\Lambda]}. 
\end{equation}

Proposition 2.9 of [BP] applied to $G' = G$ and $\Gamma' = \Gamma$ implies that any prime divisor of $[\Gamma :\Lambda]$ divides $4$, 
hence $[\Gamma:\Lambda]$ is a power of 2. Now if $\mu(\cG/\Gamma)$ is a submultiple of $1$ (i.e., it is the reciprocal of an integer) then we have the following:  

\begin{prop}
The  numerators of $\cR$ and  $\cR\prod_{v\in\cT}e'(P_v)$ are powers of $2$.   
\end{prop}


\begin{center}
{\bf 3. Discriminant bounds}
\end{center}
\vskip3mm

\ni{\bf 3.1}  We list in this section the basic estimates to be used later.
In the notation of [PY3], {\bf 2.1}, the Haar measure is normalized so that
$|\chi(\Gamma)|=\chi(X_u)\mu(\cG/\Gamma),$ where $\cG=\prod_{i=1}^rG(k_{v_i})$, $X_u$ is the 
compact dual of the symmetric space of $\cG$, and $\mu(\cG/\Gamma)$ is a submultiple of $1/{4^r}$.

We are interested in $\Gamma$ satisfying $\chi(\Gamma)\leqslant 1$, where
$\chi(\oGamma)=\chi(\Gamma)/4^r$.  We
derive from bound (4) that
\begin{equation} 1\geqslant 4^r\mu(\cG/\Gamma)\geqslant \frac{\mu(\cG/\Lambda)}{2^{d+2\#\cT}h_{\ell,4}}.
\end{equation}
From {8.1} of [PY3], we conclude that $\cE>4^{\#\cT}$ and hence we obtain the following from (1), 
\begin{equation}
1 > (D_k D_\ell)^{5/2}\big(\frac3{2^8\pi^9}\big)^d\cdot\frac1{h_{\ell,4}}.
\end{equation}

\ms
\ni{\bf 3.2} 
According to Proposition 3(ii) of [PY3], $d\leqslant 2.$ The following bound for 
$D_k^{1/d}$ is obtained from bound (57) of [PY3] for $n =3$.
\begin{equation}
D_k^{1/d} <
f_1(d,h_{\ell,4}):=\big[ \big(\frac{2^8\pi^9}3\big)^d
\cdot h_{\ell,4}
\big]^{2/15d}.
\end{equation}

We have the following bounds which are respectively  the bounds (64) and (65) of [PY3] for $n=3$: 

\begin{equation}
D_{\ell}^{1/2d}< {\mathfrak q}_1(d,D_k,h_{\ell,4})\ \ \ \ \ \ \ \ \ \ \ \ \end{equation} $$:=\big[\frac{h_{\ell,4}}{D_k^{5/2}} \cdot \big(
\frac{2^8\pi^9}3
\big)^d \big]^{1/5d}.$$

\begin{equation} 
D_{\ell}^{1/2d}< {\mathfrak q}_2(d,D_k,R_{\ell}/w_{\ell},\delta) \ \ \ \end{equation}
$$:= \big[\frac{\delta(1+\delta)}{(R_\ell/w_\ell)D_k^{5/2}}
\cdot \big( \frac{\Gamma(1+\delta)\zeta(1+\delta)^2}{(2\pi)^{1+\delta}}\cdot
\frac{2^8\pi^9}3\big)^d \big]^{1/{d(4-\delta)}}.$$
In the above, $R_\ell$ is the regulator of $\ell$ and $w_\ell$ is the order of the finite group of roots of unity contained in $\ell$.

Similarly, we have the following bounds for the relative discriminant $D_{\ell}/D_k^2$, obtained  from bounds (61) and (62) of  [PY3].

\begin{equation}
D_{\ell}/D_k^2< {\mathfrak p}_1(d,D_k,h_{\ell,4}):=\Big[h_{\ell,4}\cdot \big( \frac{2^8\pi^9}3
\big)^d D_k^{-15/2}\Big]^{2/5}.
\end{equation}

\begin{eqnarray} 
\ \ \ \ \ \ \ \ \ \ \ \ D_{\ell}/D_k^2&<& {\mathfrak p}_2(d,D_k,R_{\ell}/w_{\ell},\delta)\\ 
\nonumber&:=& \Big[\frac{\delta(1+\delta)}{(R_\ell/w_\ell)D_k^{(13-2\delta)/2}}
 \big( \frac{\Gamma(1+\delta)\zeta(1+\delta)^2}{(2\pi)^{1+\delta}}
\frac{2^8\pi^9}3\big)^d \Big]^{2/{(4-\delta)}}.
\end{eqnarray}


\bs

\begin{center}
{\bf 4. Limiting the possible pairs $(k,\ell)$}
\end{center}

\vskip3mm

 We are going to limit possibilities for the number fields $(k,\ell)$ involved in the description of  $G$.  As recalled in 3.2, $d\leqslant 2.$

\ni{\bf 4.1} As in [PY3], we are going to use extensively the lower bounds for root discriminants of number fields given in [M], which 
in turn can be traced to earlier work of Odlyzko and Diaz y Diaz.  In brief, Martinet, and earlier Odlyzko, gave increasing functions of $n$ which provide a lower bound for the root discriminant of all totally real (resp.\,totally complex) number fields of degree $n$, see ([O]; and [M], \S1.4).  
The fact that the bounds for root discriminant  provided in [O] and [M] for totally real (resp.,\,totally complex) number fields of degree $n$, increase with $n$,  has been used implicitly throughout [PY3] and will also be used here. 

We study first the case $d=2.$  
Suppose that $D_k\geqslant 33.$  In this case, $R_\ell/w_\ell$ is bounded from below by $1/8$, see [Fr, Theorem B$'$].  It follows from bound (11) that
$D_\ell^{1/4}\leqslant {\mathfrak q}_2(2,33,1/8,0.6)<18.93.$   The following argument involving Hilbert class fields will
be used repeatedly in the following.   Denote by $M_c(n)$ the smallest root discriminant among all totally complex number fields
of degree $n$.  The Hilbert class field of $\ell$ is a totally complex number field of 
degree $h_\ell$ over $\ell$, hence is of degree $4h_\ell$ over $\bQ$, with root discriminant the same as $D_\ell^{1/4}<18.93.$
On the other hand, 
it follows from Table IV of [M] that
$M_c(260)>18.98$ (hence, in fact, $M_c(n)>18.98$ for all $n\geqslant 260$).  So we conclude that $4h_\ell< 260$.
Therefore, $h_\ell\leqslant \lfloor 259/4\rfloor=64,$ where $\lfloor x\rfloor$ denotes the integral part of $x$. So $h_{\ell, 4}\leqslant 64$. 
It follows from bound (10)
that $D_\ell^{1/4}\leqslant {\mathfrak q}_1(2,33,64)<12.09.$  We iterate the above argument using Hilbert class fields.
From [M] again, $M_c(34)>12.27.$  Hence $h_\ell\leqslant \lfloor34/4\rfloor=8$. 
It follows again from bound (10)
 that $D_\ell^{1/4}\leqslant {\mathfrak q}_1(2,33,8)<9.82.$
From [M], $M_c(22)>10.25.$ 
 Hence $h_\ell\leqslant \lfloor22/4\rfloor=5$.  This implies that
$h_{\ell,4}\leqslant 4.$
Hence $D_\ell^{1/4}\leqslant {\mathfrak q}_1(2,33,4)<9.16.$

We also have the bounds 
$D_k^{1/2}\leqslant f_1(2,4)<7.84$, and 
$D_\ell/D_k^2\leqslant \lfloor{\mathfrak p}_1(2,33,4)\rfloor=6.$
Hence we have the following constraints when $D_k\geqslant 33.$
$$D_k\leqslant\lfloor7.84^2\rfloor=61, \ D_\ell \leqslant\lfloor9.16^4\rfloor=7040,\ D_\ell/D_k^2\leqslant 6.$$
As we saw in [PY1], \S8, the above constraints imply that  the possible pairs of number fields are the ones listed as $\cC_{21}-\cC_{26}$ in \S8 of [PY1].

\bs
\ni{\bf 4.2} Suppose now that $D_k<33$.  We know that the discriminant of a real quadratic number field with discriminant $<33$ is one of the following: 
$$5,\ 8,\ 12,\ 13,\ 17,\ 21,\ 24,\ 28,\ 29.$$
 We are going to handle each of these  real quadratic fields separately.

For $D_k=29$, we use the known estimates of
regulator $R_k$ provided in [C].  We know that the group of roots of unity in $\ell$ is a cyclic group of even order denoted by $m$. Let $\zeta_m$ be a primitive 
$m$-th root of unity. As the degree of the cyclotomic field $\bQ(\zeta_m)$ is $\phi(m)$, where $\phi$ is the Euler function, we know that $\phi(m)$ is a divisor of 
$2d = 4$.
Observe that the Euler function
$\phi(m)$ can take the value less than or equal to $4$ only for the following values of $m$.

$$\begin{array}{ccccccc}
m:&2&4&6&8&10&12\\
\phi(m):&1&2&2&4&4&4
\end{array}$$

Hence either $w_\ell=4$ or $w_\ell\leqslant 2$.  In the first case, as $\ell$ contains $\bQ(\zeta_m)$ and both have degree $4$, 
we conclude that $\ell=\bQ(\zeta_m),$
where $m=8,10$ or $12$.   
In the second case with $w_\ell\leqslant2,$ we know from a standard fact that 
$R_\ell/w_\ell\leqslant R_k/2,$  (cf. [W], Proposition 4.16, page 42).  Since $D_k=29$, we quote from [C], Table B.2, page 515, that $R_k=1.647$.  
Hence from estimates in (11),
\begin{equation}
D_\ell^{1/4}\leqslant {\mathfrak q}_2(2,29,1.647/2,0.65)<15.02.
\end{equation} 
Clearly the above root discriminant bound is satisfied by $\ell=\bQ(\zeta_m)$ for $m=8, 10$ or $12$ as well, cf. [W], Proposition 2.7, page 12.  
Hence the bound in (14) holds for both cases.

Since $M_c(68)>15.14$ from [M], we infer using Hilbert class field argument as above that
$h_\ell\leqslant \lfloor67/4\rfloor=16.$
This implies that $D_\ell^{1/4}\leqslant {\mathfrak q}_1(2,29,16)<11.45.$
Since $M_c(30)>11.7$ from [M],  we know that $h_\ell\leqslant \lfloor29/4\rfloor=7.$  Hence $h_{\ell,4}\leqslant 4.$
This implies that $D_\ell^{1/4}\leqslant {\mathfrak q}_1(2,29,4)<9.46.$  
 Now since $\lfloor{\mathfrak p}_1(2,29,4)\rfloor=9,$ from bound (12) we see that $D_{\ell}/D_k^2\leqslant 9$, and so 
$D_\ell\leqslant 29^2\cdot9\leqslant 7569.$  In conclusion, we have
$$D_k=29, \ D_\ell\leqslant7569,  \ D_\ell/D_k^2\leqslant 9.$$

\bs
\ni{\bf 4.3}  Employing similar arguments,
 we can handle the case of
$D_k=8, 12, 13, 17$ using the information on $R_k$ from Table B.2 of [C], on the class number $h_\ell$ from the tables in [1], and lower bounds of root discriminant from 
the tables of [M].
The explicit lower bound for $R_k$
is listed in the second
row of the table below.  

$$\begin{array}{|c|c|c|c|c|c|c|c|c|}
\hline
D_k&5&8&12&13&17&21&24&28\\ \hline
R_k>R_0=&0.4811&0.8813&1.317&1.194&2.094&1.566&2.291&2.768\\ \hline
\delta&0.492&0.534&0.573&0.576&0.615&0.623&0.649&0.672\\ \hline
{\mathfrak q}_2(2,D_k,R_0/2,\delta)<&34.1&26.4&21.5&21.2&17.7&17.1&15.4&14.1\\ \hline
{\mathfrak q}_1(2,D_k,64)<&19.37&17.3&15.6&15.3&14.3&13.6&13.1&12.6\\ \hline
h_{\ell,4}\leqslant&8&16&8&8&4&4&4&4\\ \hline
D_\ell\leqslant&61175&50458&25493&23532&13638&11040&9660&8280\\ \hline
D_\ell/D_k^2\leqslant&2447&788&177&139&47&25&16&10 \\ \hline
\end{array}
$$

\vskip1mm

Let us explain how the above table is obtained.  We first consider the cases of $D_k\geqslant 8.$
The fourth row comes from direct computation with choice of
$\delta$ given in the third row.  
Except for the
case of $D_k=5,$ the values on the fourth row satisfy
${\mathfrak q}_2(2,D_k,R_0/2,\delta)^4<10^6.$  Hence the value of $D_\ell$
lies in the list in [1].  From [1], we
check that $h_\ell\leqslant 70,$  which implies that $h_{\ell,4}\leqslant 64.$ 
This is the beginning of an argument that uses 
 Hilbert Class Fields and the tables of [M] as in the earlier cases of $D_k\geqslant33$ and 
$D_k=29$ discussed above.

In particular, this argument allows us to conclude
the upper bound $h_{\ell,4}\leqslant 4$ for $D_k\geqslant 17$ listed in the fifth row, note that here
$n=3.$  Let us carry out the procedure in detail for $D_k=8, 12, 13$ and $17.$

 For $D_k=8,$ $D_\ell^{1/2d}<{\mathfrak q}_1(2,8,64)<17.22.$ Since $M_c(130)>17.28$ from [M],  
we know that $h_\ell\leqslant \lfloor129/4\rfloor=32.$  
Hence $h_{\ell,4}\leqslant 32.$
This implies that $D_\ell^{1/2d}\leqslant {\mathfrak q}_1(2,8,32)<16.064,$  
Since $M_c(88)>16.066$ from [M],  
we know that $h_\ell\leqslant \lfloor87/4\rfloor=21.$  
Hence $h_{\ell,4}\leqslant 16.$
This implies that $D_\ell^{1/2d}\leqslant {\mathfrak q}_1(2,8,16)<14.99,$  
and $D_\ell\leqslant\lfloor{\mathfrak q}_1(2,8,16)^4\rfloor=50458.$
Furthermore, 
$D_\ell/D_k^2\leqslant \lfloor{\mathfrak p}_1(2,8,16)\rfloor=788.$  

For $D_k=12,$ $D_\ell^{1/2d}<{\mathfrak q}_1(2,12,64)<15.56.$ Since $M_c(80)>15.7$ from [M],  
we know that $h_\ell\leqslant \lfloor79/4\rfloor=19.$  Hence $h_{\ell,4}\leqslant 16.$
This implies that $D_\ell^{1/2d}\leqslant {\mathfrak q}_1(2,12,16)<13.55.$  
Since $M_c(46)>13.59$ from [M],  
we know that $h_\ell\leqslant \lfloor45/4\rfloor=11.$    Hence $h_{\ell,4}\leqslant 8.$
This implies that $D_\ell^{1/2d}\leqslant {\mathfrak q}_1(2,12,8)<12.6359$  
 and $D_\ell\leqslant\lfloor{\mathfrak q}_1(2,12,8)^4\rfloor<\lfloor12.6359^4\rfloor=
 25493.$
Furthermore, $D_\ell/D_k^2\leqslant \lfloor{\mathfrak p}_1(2,12,8)\rfloor=177.$

Since the expression ${\mathfrak q}_1(2,D_k,64)$ decreases  as $D_k$  increases,  we conclude that 
for $D_k\geqslant 12$, the argument of
the last paragraph implies that $D_\ell^{1/2d}\leqslant{\mathfrak q}_1(2,D_k,16)< {\mathfrak q}_1(2,12,16)<13.55,$
which in turn leads to $h_{\ell,4}\leqslant 8$ as well. 

In particular, for $D_k=13$, $h_{\ell,4}\leqslant 8$.

For  $D_k\geqslant 17$,  $D_\ell^{1/2d}\leqslant{\mathfrak q}_1(2,D_k,8)\leqslant {\mathfrak q}_1(2,17,8)<11.5823.$
Since $M_c(30)>11.70$ from [M],  
we know that $h_\ell\leqslant \lfloor29/4\rfloor=7.$   It follows  that $h_{\ell,4}\leqslant 4.$  The other items in the 
above table are computed similarly.

\vskip1mm

For $D_k=5,$ ($k = \bQ({\sqrt{5}})$) we simply use the bound
$D_\ell\leqslant {\mathfrak q}_2(2,5,0.4811/2,0.5)^4<34.1^4<1.36\times 10^6.$  
At the request of the authors, Malle provide us with a complete list of all totally complex number fields of degree 4 containing $\bQ({\sqrt{5}})$ satisfying this  discriminant bound. 
There are $2556$ such number fields and for each $\ell$ in this list, $h_\ell\leqslant 65.$  Hence $h_{\ell,4}\leqslant 64.$
It follows that 
$D_\ell\leqslant {\mathfrak q}_1(2,5,64)^{4}\leqslant 140565.$  From the list of Malle, it follows that there are $276$ totally complex number fields
satisfying discriminant $D_\ell\leqslant 140565$.  Moreover, $h_\ell\leqslant 18$ and hence $h_{\ell,4}\leqslant 16,$
from which we conclude that $D_\ell\leqslant \lfloor{\mathfrak q}_1(2,5,16)^{4}\rfloor\leqslant 80733.$
Checking in Malle's list of number fields again, there are $164$ number fields $\ell$ with discriminant $D_\ell\leqslant 80733$, we find that $h_\ell\leqslant 12$ for them, and hence $h_{\ell,4}\leqslant 8.$
We conclude that $D_\ell/D_k^2\leqslant \lfloor{\mathfrak p}_1(2,5,8)\rfloor\leqslant 2447.$  Hence 
$D_\ell\leqslant 61175$.  There are altogether $121$ such number fields $\ell$.



\bs
\ni{\bf 4.4} We summarize the results from the previous subsections in the following.  

\begin{prop} The followings are all possible pairs $(k,\ell)$ which may give 
an arithmetic fake compact hermitian symmetric space of type $A_3$.\\
(a)  $[k:\bQ]=2$, and $(k,\ell)$ is one of the pairs $\cC_{21}-\cC_{26}$ described in \S8 of [PY1], 
or $k =\bQ(\sqrt{\alpha})$ for $\alpha$ as in the table below, and the discriminant $D_{\ell}$ of $\ell$ satisfies the following bound:

$$\begin{array}{|c|c|c|c|c|c|c|c|c|c|}
\hline
\alpha&5&2&3&13&17&21&6&7&29\\  \hline
D_\ell\leqslant &61175&50458&25493&23532&13638&11040&9660&8280&7569 \\ \hline
N &121&50&13&12&4&3&2&2&1\\ \hline
\end{array}$$
where $N$ is the number of number fields satisfying the discriminant bounds.  

\vskip2mm
\ni(b) $k=\bQ,$ and $\ell=\bQ(\sqrt{-a})$:  \ \ 
$D_\ell\leqslant1363.$
\end{prop}

Assertion (b) is from {9.3} of [PY3].  We note that in case (b) there are altogether $434$ possible number fields $\ell$.  In (a), there are 
altogether $208$ pairs of number fields in the table, apart from the $6$ pairs of number fields $\cC_{21}-\cC_{26}$.  The authors are grateful
to Gunter Malle for providing the list of number fields satisfying the constraints in the table above.  The list of number fields can be found in the weblink
[2] provided by Malle.

\vskip4mm

\begin{center}
{\bf 5.  Zeta and $L$-values}
\end{center}
\vskip3mm

 In this and the next section, we are going to use the value of $\cR$ given by (5) and the explicit values of $\zeta$- and $L$-functions to restrict the possible number fields involved.  We will 
use the fact that the numerators of $\cR$ and $\cR\prod_{v\in\cT}e'(P_v)$ are powers of $2$ as given by Proposition 1.  

\ms
\ni{\bf 5.1} Consider $\cC_i,$ for $i= 21$--$26$ the pairs of number fields $(k,\ell)$ introduced in [PY1], \S8.  Denote by $\zeta_n$ a primitive $n$-th 
root of unity.  Then 
$$\cC_{21}=(\bQ({\sqrt{33}}), \bQ({\sqrt{33}},\zeta_5)), \cC_{22}=(\bQ({\sqrt{11}}), \bQ({\sqrt{11}},\zeta_4)),$$ 
$$\cC_{23}=(\bQ({\sqrt{14}}), \bQ({\sqrt{14}},{\sqrt{-7}})),
\cC_{24}=(\bQ({\sqrt{57}}), 
\bQ({\sqrt{57}},\zeta_3)),$$ 
$$\cC_{25}=(\bQ({\sqrt{15}}), 
\bQ({\sqrt{15}},\sqrt{-5})), \  
 \cC_{26}=(\bQ({\sqrt
{15}}), \bQ({\sqrt{15}},\zeta_4)).$$
From computations using Magma, we have
the following table of values for $\zeta_k(-1),$ $L_{\ell/k}(-2),$  $\zeta_k(-2)$ \ \ and \ \ $\cR$.
$$\begin{array}{|c|c|l|l|l|}
\hline
(k,\ell)&\zeta_k(-1)&L_{\ell|k}(-2)&\zeta_k(-2)&\cR\\  \hline\hline
\cC_{21}&1&4/3&141/10&47/160\\  \hline
\cC_{22}&7/6&3&2153/60&15071/7680\\  \hline
\cC_{23}&5/3&48/7&2503/30&2503/168\\  \hline
\cC_{24}&7/3&44/9&2867/30&220759/12960\\  \hline
\cC_{25}&2&20/3&537/5&179/8\\  \hline
\cC_{26}&2&8&537/5&537/20\\  \hline
\end{array}$$

By looking at the last column of the above table, we see that the cases $\cC_{21}-\cC_{26}$ can be ruled out since the numerators
are not power of $2$.

\ms
\ni{\bf 5.2} 
 We compute the value of $\cR$ for each of the candidate  pairs  $(k,\ell)$ provided by Malle in Proposition 2(a).  Only four of them have numerator a power of $2.$   These are listed
below.  Note that for $k=\bQ(\sqrt{5})$ the zeta values are 
$\zeta_k(-1)=1/30$ and $\zeta_k(-3)=1/60.$

$$\begin{array}{|c|c|l|}
\hline
(k,\ell)&L_{\ell|k}(-2)&\cR\\  \hline\hline
\cC_1&4/5&1/144000=1/(2^7\cdot3^2\cdot5^3)\\  \hline
\cC_2&32/9&1/32400=1/(2^4\cdot3^4\cdot5^2)\\  \hline
\cC_3&15&1/7680=1/(2^9\cdot3\cdot 5)\\  \hline
\mathscr{F}_1&2400&1/48=1/(2^4\cdot3) \\  \hline
\end{array}$$

Here  $\cC_i,$ $i=1,2,3$ are pairs as  in [PY1], \S8,   and ${\mathscr{F}}_1= (\bQ({\sqrt 5}),\bQ({\sqrt 5},\sqrt{-11}))$, 
with $D_k=5$ and $D_\ell=3025$. Recall that $\cC_1=(\bQ({\sqrt 5}), \bQ(\zeta_5))$; $\cC_2= (\bQ({\sqrt 5}), \bQ({\sqrt 5}, \zeta_3 ))$,
and  $\cC_3= (\bQ({\sqrt 5}), \bQ({\sqrt 5}, \zeta_4 ))$.
Note that the class number of $\ell$ in $\cC_1$, $\cC_2$ and $\cC_3$ is $1$, while the class number of $\ell\ (=\bQ({\sqrt{5}}, {\sqrt{-11}}))$ in $\mathscr{F}_1$ is $2.$


\vskip2mm

\ni{\bf 5.3}  Consider now the case of $k=\bQ.$  
In this case, $\zeta_k$ is just the regular Riemann Zeta function.

Computing the values of $L_{\ell/{\bQ}}(-2)$ for $\ell$ 
with $D_\ell\leqslant 1363,$ we see that 
the only candidates for the cases of $k=\bQ,\ \ell=\bQ(\sqrt{-a})$ are the following ten for which the numerator of $\cR$ is a power of $2$. In this table,  $h$ is the class number of $\ell=\bQ(\sqrt{-a})$.
$$
\begin{array}{|c|c|l|}
\hline 
a&h&\cR\\ \hline\hline
1&1&1/23040=1/(2^9\cdot3^2\cdot5)\\ \hline
2&1&1/3840=1/(2^8\cdot3\cdot5)\\  \hline
3&1&1/51840=1/(2^7\cdot3^4\cdot5)\\  \hline
5&2&1/384=1/(2^7\cdot3)\\  \hline
7&1&1/5040=1/(2^4\cdot3^2\cdot5\cdot7)\\  \hline
11&1&1/1920=1/(2^7\cdot3\cdot5)\\  \hline
15&2&1/720=1/(2^4\cdot3^2\cdot5)\\  \hline
23&3&1/240=1/(2^4\cdot3\cdot5)\\  \hline
31&3&1/120=1/(2^3\cdot3\cdot5)\\  \hline
47&5&1/40=1/(2^3\cdot5)  \\ \hline
\end{array}
$$

\bs
\begin{center}
{\bf 6. Contributions by $e'(P_v)$}
\end{center}

\bs
\ni{\bf 6.1}  We are going to make use of the factor $e'(P_v)$ in $\mu(\mathcal{G}/\Lambda)= \cR\prod_{\cT}e'(P_v)$ to eliminate some more cases.

We will consider only the cases where the division algebra $\cD\ne \ell$. Then $\cT_0$ is nonempty and every place in $\cT_0$ is nonarchimedean and splits in $\ell$ (see 1.1)
\vskip1mm

Here is the list of  rational primes  
$p\leqslant 71$ which split in $\ell$.  {\it Note that for the pairs $\cC_1$, $\cC_2$, $\cC_3$ and $\mathscr{F}_1$, we only list those  $p\leqslant 71$, which are restrictions to $\bQ$ of nonarchimedean places of  $k$ which split in $\ell$.} For any nonarchimedean place $v$ of $k$ which splits in $\ell$ and lies over a rational prime larger than $71$, the value of $\cR e'(P_v)\gg 1$ for all the 14 cases presently under consideration (these are the cases for which $\cR$ has been listed in 5.2 and 5.3).

\ms
$$
\begin{array}{|c|l|}
\hline
(k,\ell) & \mbox{primes $\leqslant 71$ which split in $\ell$}\\ \hline\hline
(\bQ,\bQ(\sqrt{-1}))&5,13,17,29,37,41,53,61\\ \cline{2-2}
(\bQ,\bQ(\sqrt{-2}))&3,11,17,19,41,43,59,67\\ \cline{2-2}
(\bQ,\bQ(\sqrt{-3}))&7,13,19,31,37,43,61,67\\ \cline{2-2}
(\bQ,\bQ(\sqrt{-5}))&3,7,23,29,41,43,47,61,67\\ \cline{2-2}
(\bQ,\bQ(\sqrt{-7}))&2,11,23,29,37,43,53,67,71\\ \cline{2-2}
(\bQ,\bQ(\sqrt{-11}))&3,5,23,31,37,47,53,59,67,71\\ \cline{2-2}
(\bQ,\bQ(\sqrt{-15}))&2,17,19,23,31,47,53,61\\ \cline{2-2}
(\bQ,\bQ(\sqrt{-23}))&2,3,13,29,31,41,47,59,71 \\ \cline{2-2}
(\bQ,\bQ(\sqrt{-31}))&2,5,7,19,41,47,59,67,71\\ \cline{2-2}
(\bQ,\bQ(\sqrt{-47}))&2,3,7,17,37,53,59,61,71\\ \hline
\cC_1=(\bQ({\sqrt 5}), \bQ(\zeta_5))&11,31,41,61,71\\ \cline{2-2}
\cC_2=(\bQ({\sqrt 5}), \bQ({\sqrt 5}, \zeta_3 ))&7,13,17,19,23,31,37,43,53,61,67\\ \cline{2-2}
\cC_3=(\bQ({\sqrt 5}), \bQ({\sqrt 5}, \zeta_4 ))&3,5,7,13,17,23,29,37,41,43,47,53,61,67\\ \cline{2-2}
{\mathscr{F}}_1=(\bQ({\sqrt 5}), \bQ({\sqrt 5},\sqrt{-11}))&3,5,7,13,17,23,31,37,43,47,53,59,67,71\\ \hline
\end{array}
$$

The primes $\leqslant 71$ which split in $\bQ({\sqrt{5}})$ are 11, 19,  29, 31, 41, 59, 61 and 71.

\ms
Note that in each of the cases above, $\ell$ is a Galois extension of $\bQ$.
\vskip2mm

\ni{\bf 6.2} The values of $e'(P_v)$ for $v\in \cT_0$ are given in Case 2 of {\S2.3}.  We see that $e'(P_v)$ is an integral multiple of $(q_v-1)(q^3_v-1).$  By direct computation,
we conclude that except for the pair $(\bQ, \bQ(\sqrt{-7}))$, for every other pair in the above table, for any possible choice of a nonempty 
$\cT_0$ either $ \cR\prod_{v\in\cT_0}e'(P_v)\gg 1$ or its numerator is not a power of $2$.   We also observe that for $(\bQ,\bQ(\sqrt{-7}))$, 
the only possibility for  $\cT_0$ is $\cT_0=\{2\}$. 
\vskip1mm

We know from 9.6 of [PY3] that if $k =\bQ$, then $\cD \ne \ell$. Thus the only possible pair with $k =\bQ$ is $(\bQ, \bQ(\sqrt{-7}))$ and in this case $\cT_0=\{ 2\}$ and  $\sqrt{[\cD:\ell]} =2$ or $4$.

\bs
We summarize the above results as follows.

\begin{theo} The following are all possible pairs $(k,\ell)$ which may give rise to 
a fake compact hermitian symmetric space of type $A_3$.\\
(a) $(k,\ell)\in\{\cC_1, \cC_2, \cC_3, \mathscr{F}_1\}$ and $\cD = \ell$.\\
(b) $(k,\ell)=(\bQ, \bQ(\sqrt{-7}))$, ${\sqrt{[\cD:\ell]}}=2$ or $4$ and $\cT_0 = \{2\}$.
\end{theo}

\bs
\begin{center}
 {\bf 7. Potential examples for $k=\bQ$.}
\end{center}

\ms 
\bs
\ni{\bf 7.1} We consider the case of $(k,\ell)=(\bQ, \bQ(\sqrt{-7}))$. If ${\sqrt{[\cD:\ell]}}=2$, i.e., $\cD$ is a quaternion division algebra, then $\cT_0$ determines it uniquely.  Using the values of $e'(P_2)$ given in Cases 2b(i) and 2b(ii) in \S2.3 we see that either $P_2$ can be a maximal parahoric subgroup  in which case $e'(P_2) = 7$ and hence $\cR e'(P_2) = 7/5040=1/720$, or $P_2$ can be an Iwahori subgroup in which case $e'(P_2) = 35$ and $\cR e'(P_v)= 35/5040 = 1/144$.     
\vskip1mm

Now since $G$ is anisotropic over $k$,  
there exist at least two primes $q$ which do not split in $\bQ(\sqrt{-7})$ and $\bQ_q$-rank of $G$ is $1$; in fact the number of such primes is even (and at least 2).  
Let
\begin{eqnarray*}
\psi_1(q)&=&(q^2-1),\\
\psi_2(q)&=&(q-1)(q^2+1).
\end{eqnarray*}
From the values of $e'(P_v)$ given in Cases 3b, 4b of \S2.3,  if $q$ is ramified in $\ell=\bQ(\sqrt{-7})$ (i.e., $q =7$), $e'(P_7)$ is an integral multiple of either $\psi_1(7)= 48$ or it is $\psi_2(7)= 300$, and if it is unramified in $\ell$, it is an integral multiple of $\psi_2(q)$. 
\vskip1mm

For small primes $q$ which do not split in $\bQ(\sqrt{-7})$, we have the following values of  $\psi_2(q).$  Clearly, both
these functions are increasing in $q.$  So we find that it suffices for us to consider $\psi_2(q)$ only
for $q=3,\ 5, \ 7$ and $13.$

$$\begin{array}{|c|c|}
\hline
\psi_2(3)&2^2\cdot 5\\ \hline
\psi_2(5)&2^3\cdot 13\\ \hline
\psi_2(7)&2^2\cdot 3\cdot 5^2\\ \hline
\psi_2(13)&2^3\cdot 3\cdot 5\cdot 17\\ \hline
\end{array}$$
\vskip2mm

By computing $\cR e'(P_2)\prod_{p\in \cT-\{2\}}e'(P_p)$ and checking whether  its numerator is a power of 2 (recall from the above that $\cT$ contains at least two primes $q$ 
which do not split in $\bQ(\sqrt{-7})$ and $\bQ_q$-rank of $G$ is $1$), we see that $\cT$ must equal $\{2, 3, 7\}$ and moreover, for each $p\in \cT$, $P_p$ is a maximal parahoric subgroup of $G(\bQ_p)$ 
and $\oM_7 =R_{\fF_7/\ff_7}({\rm{SL}}_2)$. This data determines a principal arithmetic subgroup $\Lambda$ in ${\rm{SU}}(2,2)$ whose covolume is $\cR e'(P_2)e'(P_3)e'(P_7) = 4/3$. Now the question is whether its normalizer in 
${\rm{PU}}(2,2)$ contains a torsion-free subgroup of suitable index.  

\bs

\ni{\bf 7.2} Let us now consider the case where the degree of $\cD$ is $4$. Using the formulae in \S2.3, Case 2a, we see that one possibility is $\cT = \cT_0 =\{2\}$, in which case $\cR e'(P_2)=21/ 5040=1/240.$  (In this case $\oM_2=R_{\mathfrak{F}_2/\ff_2}({\rm{GL}}_1)/{\rm{GL}}_1$, determines a principal arithmetic subgroup $\Lambda'$ of covolume $1/240$ in ${\rm{SU}}(2,2)$.  
\vskip1mm

\begin{ques} Does  $\Lambda'$ contain a torsion-free subgroup of index $240$?
\end{ques}

  Let $\Lambda$  be the second congruence subgroup ${\rm{SL}}_1^{(2)}(\bQ_2\otimes_{\bQ}\cD)$
of $G(\bQ_2)={\rm{SL}}_1(\cD)$; it  is a normal subgroup of index $240$. However, this $\Lambda$ does contain elements of order 4. 
 
\vskip1mm

\ni{\bf 7.3} Another possibility (when $\cD$ is of degree 4) is that  at primes $3$ and $7$ the group $G$ is of rank 1 and the parahoric subgroup $P_p$ for $p =3, 7$ is maximal with  $\oM_7 =R_{\fF_7/\ff_7}({\rm{SL}}_2)$. In this case, $\cT =\{2, 3, 7\}$ and we obtain a principal arithmetic subgroup $\Lambda''$ whose covolume is $e'(P_3)e'(P_7)/240= 20\cdot 48/240 = 4$.

\bs
\begin{center}
{\bf 8 Potential examples for $[k,\bQ]>1$.}
\end{center}

\ms
\ni{\bf 8.1} We consider the case of $(k,\ell)=\cC_1, \cC_2, \cC_3, \mathscr{F}_1$ and $\cD=\ell$.  Then there exists a hermitian form $h$ on $\ell^4$ (defined in terms of the nontrivial automorphism of $\ell/k$) 
such that $G={\rm{SU}}(h)$.
 
As $k = \bQ({\sqrt{5}})$ is a field of degree $2$,  $r \leqslant 2$.  Here we claim that actually $r=1$. To prove this claim, we assume on the contrary that $r=2$.  
Then $G$ is isotropic at both the real places of $k = \bQ({\sqrt{5}})$. On the other hand,  for a place $v$ of $k$, the group 
$G$ is isomorphic to the split group ${\rm{SL}}_4$ over $k_v$ if $v$ splits in $\ell$, and if $v$ does not split in $\ell$, then $G$ is $k_v$-isomorphic 
to the special unitary group of a hermitian form on $\ell_v^4$, where  $\ell_v = \ell\otimes_k k_v$.   So we see, for example, from Proposition 7.2 of [PR] that 
$G$ is isotropic over $k$ and hence we get a non-compact locally Hermitian symmetric space, 
contradicting our initial assumption.

\ms
\ni{\bf 8.2}  For the four pairs of fields $\cC_1$, $\cC_2$, $\cC_3$ and $\mathscr{F}_1$, we will now list the parahoric subgroups involved in the description of the principal arithmetic subgroup $\Lambda$. Let us assume that $(k,\ell)$ is one of the four pairs under consideration. Let us denote by $\cT_1$ the set of nonarchimedean places $v$ of $k$ which do not split in $k$ and $G$ is of $k_v$-rank 1. It is known that the cardinality of $\cT_1$ is even, $\cT_1$ can be empty. We are going to show that in fact $\cT_1$ is empty and thus the Witt index of $h$ is $2$ at every nonarchimedean place of $k$ which does not split in $\ell$.  To see this let us assume that $\cT_1$ is nonempty, and  let $v'$ and $v''$ be two places belonging to it. By computing $\cR e'(P_{v'})e'(P_{v''})$ using the values listed in Cases 3b and 4b in \S2.3 we see that if either of these places lies over a rational prime $>29$, then $\cR e'(P_{v'})e'(P_{v''})\gg 1$. On the other hand, if $v'$ and $v''$ lie over primes  $\leqslant 29$, then the numerator of $\cR e'(P_{v'})e'(P_{v''})$  is not a power of $2$. From these observations we conclude that $\cT_1$ is empty.

Also, by computing $\cR e'(P_v)$ using the values of $e'(P_v)$ given in \S2.3, Cases 1, 3a and 4a, we see that $P_v$  has to be hyperspecial if either $v$ splits in $\ell$ or is inert, and if $v$ ramifies in $\ell$, then $P_v$ is a special maximal parahoric subgroup. So $e'(P_v) =1$ for all nonarchimedean $v$. This implies that the covolume of $\Lambda$ is precisely $\cR$.

\bs
\begin{ques}
Does there exist torsion-free subgroup of right index, so that its covolume is $1$, in the normalizer $\Gamma$ of $\Lambda$  in ${\rm{SU}}(2,2)$ for \,$(k,\ell)= \cC_1, \, \cC_2, \, \cC_3,\,  \mathscr{F}_1$ 
 respectively for Hermitian form case?   We can also work with the image of $\Gamma$ in ${\rm{PU}}(2,2)$. 
 \end{ques}

\bs
The values of $\cR$ are given in {5.2}.  To rule out examples in the case of $\fF_1$, it suffices for us to find in the principal arithmetic group $\Lambda$
a torsion element of order $5$, $7$, $3^2$, or contain
a factor not dividing $2^a\cdot 3$ for any positive integer $a$.

\bs
\centerline{\bf References}

\vskip4mm

\ni[BP] A.\,Borel and G.\,Prasad, {\it Finiteness theorems for discrete subgroups of bounded covolume in semisimple groups.} Publ.\,Math.\,IHES No.\,{\bf 69}(1989), 119--171.
\vskip1.5mm

\ni [C]  H.\,Cohen, {\it A course in computational algebraic number theory.} Graduate Texts in Mathematics, 138. Springer-Verlag, Berlin, 1993.
\vskip1.5mm

\ni [CS] Cartwright, D., Steger, T., Enumeration of the $50$ fake projective planes, C. R. Acad. Sci. Paris, Ser. 1, 348 (2010), 11-13,
\vskip1.5mm

\ni [M] J.\,Martinet,  {\it Petits discriminants des corps de nombres.} Number theory days (Exeter, 1980), 151--193, London Math. Soc. Lecture Note Ser., 56, Cambridge Univ. Press, Cambridge-New York, 1982.
\vskip1.5mm

\ni [O] A.\,M.\,Odlyzko,  {\it Discriminant bounds}, unpublished, available
from:\\
http://www.dtc.umn.edu/$\sim$odlyzko/unpublished/index.html.
\vskip1.5mm
 
\ni[P] G.\,Prasad, {\it Volumes of $S$-arithmetic quotients of semi-simple groups.} Publ. Math. IHES No. {\bf 69}(1989), 91-117.
\vskip1.5mm

\ni [PR], G.\,Prasad and A.\,S.\,Rapinchuk, {\it Weakly commensurable arithmetic groups and isospectral locally symmetric spaces}, Publ.\,Math.\,IHES No.\,{\bf 109}(2009), 113-184.
\vskip1.5mm

\ni[PY1] G.\,Prasad and S-K.\,Yeung, {\it Fake projective planes}. Inv.\,Math.\,{\bf 168}(2007), 321-370. {\it Addendum}, ibid {\bf 182}(2010), 213-227. (See the corrected version posted on arXiv.) 

\vskip1.5mm
\ni
[PY2] G.\,Prasad and S-K.\,Yeung, {\it Arithmetic fake projective spaces and arithmetic fake Grassmannians}, American J.\,Math.\,{\bf 131}(2009), 379-407.

\vskip1.5mm 
\ni [PY3] G.\,Prasad and S-K.\,Yeung, {\it Nonexistence of arithmetic fake compact hermitian symmetric spaces of type other than $A_n\: (n\leqslant 4)$}, J. Math. Soc. Japan 64 (2012), no. 3, 683-731.

\vskip1.5mm 
\ni [W]  L.\,C.\,Washington, Introduction to Cyclotomic Fields, 2nd edition, Grad. Texts in Math., 83, Springer-Verlag, New York, 1997.

\vskip1.5mm 
\ni[1] The Bordeaux Database, Tables obtainable from: \\
ftp://megrez.math.u-bordeaux.fr/pub/numberfields/.

 \vskip1.5mm 
\ni[2] Malle's Database, from:\\
 http://www.mathematik.uni-kl.de/~numberfieldtables/KT\_4/download.html
\end{document}